\documentclass{amsart}
\usepackage{amssymb}
\usepackage[all]{xy}
\usepackage{hyperref}
\usepackage{a4wide}

\newtheorem{thm}{Theorem}

\newtheorem{theorem}{Theorem}
\newtheorem{lemma}[theorem]{Lemma}
\newtheorem{corollary}[theorem]{Corollary}
\newtheorem{conjecture}[theorem]{Conjecture}

\renewcommand{\phi}{\varphi}
\newcommand{\alp}{\alpha}
\newcommand{\sig}{\sigma}
\newcommand{\Del}{\Delta}
\newcommand{\bfsig}{\mbox{\boldmath$\sigma$}}
\newcommand{\bfa}{\mbox{\boldmath$a$}}
\newcommand{\Gal}{\textnormal{Gal}}
\newcommand{\bbA}{\mathbb{A}}
\newcommand{\bbQ}{\mathbb{Q}}
\newcommand{\bbZ}{\mathbb{Z}}
\newcommand{\Kgal}{{\tilde K}}
\newcommand{\Ogal}{{\tilde O}}
\newcommand{\calS}{\mathcal{S}}
\newcommand{\discr}{{\rm discr}}
\def\irr{{\rm irr}}
\def\Egag{{\bar{E}}}
\def\Fgag{{\bar{F}}}
\def\Fhat{{\hat F}}
\def\Shat{{\hat S}}
\def\Ker{{\rm Ker}}
\def\st{{\mid}}
\def\tensor{\otimes}
\def\union{\bigcup}
\def\normal{\triangleleft}
\newcommand{\isom}{\cong}
\newcommand{\cont}{\subseteq}
\newcommand{\hefresh}{\smallsetminus}

\begin{document}
\fontsize{12}{18}\selectfont
\vfuzz2pt 
\hfuzz2pt 

\title{PAC Fields over Finitely Generated Fields} 
\author{Lior Bary-Soroker \and Moshe Jarden}

\thanks{MSC-class: 12E30}

\thanks{Research supported by the Minkowski Center for Geometry at Tel Aviv
University, established by the Minerva Foundation.}

\thanks{This work constitutes a part of the Ph.D dissertation of the
first author done at Tel Aviv University under the supervision of
Prof. Dan Haran.}

\address{School of Mathematics, Tel Aviv University,\\ Ramat Aviv,
Tel Aviv 69978, Israel} \email{barylior@post.tau.ac.il and
jarden@post.tau.ac.il}
\date{\today}
%

%
\begin{abstract}
We prove the following theorem for a finitely generated field $K$:
Let $M$ be a Galois extension of $K$ which is not separably closed.
Then $M$ is not PAC over $K$.
\end{abstract}
\maketitle
\section{Introduction}
\label{intro} A central concept in Field Arithmetic is ``pseudo
algebraically closed (abbreviated {\bf PAC}) field''. If $K$ is a
countable Hilbertian field, then $K_s(\bfsig)$ is PAC for almost all
$\bfsig\in\Gal(K)^e$ \cite[Thm.~18.6.1]{FrJ}. Moreover, if $K$ is
the quotient field of a countable Hilbertian ring $R$
(e.g.~$R=\bbZ$), then $K_s(\bfsig)$ is PAC over $R$
\cite[Prop.~3.1]{JaR1}, hence also over $K$.

Here $K_s$ is a fixed separable closure of $K$ and
$\Gal(K)=\Gal(K_s/K)$ is the absolute Galois group of $K$. This
group is equipped with a Haar measure and ``almost all'' means ``for
all but a set of measure zero''. If
$\bfsig=(\sig_1,\ldots,\sig_e)\in\Gal(K)^e$, then $K_s(\bfsig)$
denotes the fixed field in $K_s$ of $\sig_1,\ldots,\sig_e$.

Recall that a field $M$ is said to be {\bf PAC} if every nonempty
absolutely irreducible variety $V$ defined over $M$ has an
$M$-rational point. One says that $M$ is {\bf PAC over} a subring
$R$ if for every absolutely irreducible variety $V$ defined over $M$
of dimension $r\ge0$ and every dominating separable rational map
$\phi\colon V\to\bbA_M^r$ there exists an $\bfa\in V(M)$ with
$\phi(\bfa)\in R^r$.

When $K$ is a number field, the stronger property of the fields
$\Kgal(\bfsig)$ (namely, being PAC over the ring of integers $O$ of
$K$) has far reaching arithmetical consequences. For example,
$\Ogal(\bfsig)$ (= the integral closure of $O$ in $\Kgal(\bfsig)$)
satisfies Rumely's local-global principle \cite[special case of
Cor.~1.9]{JaR2}: If $V$ is an absolutely irreducible variety defined
over $\Kgal(\bfsig)$ with $V(\Ogal)\ne\emptyset$, then $V$ has an
$\Ogal(\bfsig)$-rational point. Here $\Kgal$ denotes the algebraic
closure of $K$ and $\Kgal(\bfsig)$ is, as before, the fixed field of
$\sig_1,\ldots,\sig_e$ in $\Kgal$.

The article \cite{JaR1} gives several distinguished Galois
extensions of $\bbQ$ which are not PAC over any number field and
notes that no Galois extension of a number field $K$ (except
$\Kgal$) is known to be PAC over $K$. This lack of knowledge has
come to an end in \cite{Jar}, where Neukirch's characterization of
the $p$-adically closed fields among all algebraic extensions of
$\bbQ$ is used in order to prove the following theorem:

\begin{thm}\label{thm:A}
If $M$ is a Galois extension of a number field $K$ and $M$ is
not algebraically closed, then $M$ is not PAC over $K$.
\end{thm}

The goal of the present note is to generalize Theorem~\ref{thm:A} to
an arbitrary finitely generated field (over its prime field):

\begin{thm}\label{thm:B}
Let $K$ be a finitely generated field  and $M$ a Galois extension of
$K$ which is not separably closed. Then $M$ is not PAC over $K$.
\end{thm}

The proof of Theorem~\ref{thm:B} is based on Proposition 5.4 of
\cite{JaR1} which combines Faltings' theorem in characteristic $0$
and the Grauert-Manin theorem in positive characteristic. The latter
theorems are much deeper than the result of Neukirch used in the
proof of Theorem~\ref{thm:A}.

\section{Accessible extensions}
The proof of Theorem~\ref{thm:B} actually gives a stronger theorem:
No accessible extension (see definition prior to
Theorem~\ref{thm:4}) of a finitely generated field $K$ except $K_s$
is PAC over $K$. Technical tools in the proof are the ``field
crossing argument'' and ``ring covers'':

An extension $S/R$ of integral domains with an extension $F/E$ of
quotient fields is said to be a {\bf cover of rings} if $S=R[z]$ and
$\discr(\irr(z,E))\in R^\times$ \cite[Definition 6.1.3]{FrJ}. We say
that $S/R$ is a {\bf Galois cover of rings} if $S/R$ is a cover of
rings and $F/E$ is a Galois extension of fields. Every epimorphism
$\phi_0$ of $R$ onto a field $\Egag$ extends to an epimorphism
$\phi$ of $S$ onto a Galois extension $\Fgag$ of $\Egag$ and $\phi$
induces an isomorphism of the {\bf decomposition group}
$D_\phi=\{\sig\in\Gal(F/E)\st\sig(\Ker(\phi))=\Ker(\phi)\}$ onto
$\Gal(\Fgag/\Egag)$ \cite[Lemma~6.1.4]{FrJ}. In particular,
$\Gal(F/E)\isom\Gal(\Fgag/\Egag)$ if and only if
$[F:E]=[\Fgag:\Egag]$.

As in the proof of \cite[Lemma 24.1.1]{FrJ}, the field crossing
argument is the basic ingredient of the construction included in the
proof of the following lemma.

\begin{lemma}\label{lem:1}
Let $K$ be a field, $M$ an extension of $K$, $n$ a positive integer,
$N$ a Galois extension of $M$ with Galois group $A$ of order at most
$n$, and $t$ an indeterminate. Then there exist fields
$D,F_0,F,\Fhat$ as in diagram (\ref{eq1}) such that the following
holds:
\renewcommand{\labelenumi}{\textnormal{(\alph{enumi})}}
\begin{enumerate}
\item
$F_0$ is regular over $K$, $F$ and $D$ are regular over $M$, and
$\Fhat$ is regular over $N$.
\item
$FD=DN=\Fhat$.
\item
$F_0/K(t)$, $F/M(t)$, and $\Fhat/N(t)$ are Galois extensions with
Galois groups isomorphic to $S_n$.
\end{enumerate}
\begin{equation}\label{eq1}
\xymatrix{%
F_0 \ar@{-}[r] \ar@{-}[d]_{S_n}
    & F \ar@{-}[r]^A \ar@{-}[d]_{S_n}
        & \Fhat \ar@{-}[d]^{S_n} \ar@{-}[dl]|D
\\
K(t) \ar@{-}[r] \ar@{-}[d]
    & M(t) \ar@{-}[r] \ar@{-}[d]
        & N(t)\ar@{-}[d]
\\
K \ar@{-}[r]
    & M \ar@{-}[r]^A & N }
\end{equation}
\end{lemma}

\begin{proof}
The field $K(t)$ has a Galois extension $F_0$ with
Galois group $S_n$ such that $F_0$ is regular over $K$ \cite[Example
16.2.5 and Proposition 16.2.8]{FrJ}. In particular, $F_0$ is
linearly disjoint from $N$ and $M$ over $K$. Set $F=F_0M$ and
$\Fhat=FN$. By \cite[Cor.~2.6.8]{FrJ}, both $F/M$ and $\Fhat/N$ are
regular extensions. Moreover, both $F/M(t)$ and $\Fhat/N(t)$ are
Galois extensions with Galois groups isomorphic to $S_n$ and
$\Fhat/F$ is a Galois extension. We identify $\Gal(\Fhat/F)$ with
$A$ via restriction. Finally, $\Fhat/M(t)$ is a Galois extension and
$\Gal(\Fhat/M(t))=\Gal(\Fhat/F)\times\Gal(\Fhat/N(t))$.

Multiplication from the right embeds $A$ into $S_m$, where $m=|A|$.
Since $m\le n$, there exists an embedding $\alp\colon
A\to\Gal(\Fhat/N(t))$. Consider the diagonal subgroup
$\Del=\{(\sig,\alp(\sig))\in\Gal(\Fhat/M(t))\st\sig\in A\}$ of
$\Gal(\Fhat/M(t))$ and its fixed field $D$ in $\Fhat$. Then
$\Del\cap\Gal(\Fhat/F)=\Del\cap\Gal(\Fhat/N(t)=1$. By Galois theory,
$FD=DN(t)=\Fhat$, so $DN=\Fhat$. Restriction to $N$ maps
$\Gal(\Fhat/D)$ onto $\Gal(N/M)$, hence $D\cap N=M$. Since $\Fhat$
is regular over $N$, it follows that $D$ is regular over $M$.
\end{proof}

The main ingredient in the proof of Lemma~\ref{lem:3} is the
following result of Faltings' theorem in characteristic $0$ and the
Grauert-Manin theorem in positive characteristic.

\begin{lemma}[{\cite[Prop.~5.4]{JaR1}}]\label{lem:2}
Let $K$ be an infinite finitely generated field, $f\in K[T,Y]$ an
absolutely irreducible polynomial which is separable in $Y$, $g\in
K[T,Y]$ an irreducible polynomial which is separable in $Y$, and
$0\ne r\in K[T]$. Then there exist a finite purely inseparable
extension $K'$ of $K$, a nonconstant rational function $q\in K'(T)$,
and a finite subset $B$ of $K'$ such that $f(q(T),Y)$ is absolutely
irreducible, $g(q(a),Y)$ is irreducible in $K'[Y]$, and
$r(q(a))\ne0$ for each $a\in K'\hefresh B$.
\end{lemma}

\begin{lemma}\label{lem:3}
Let $K$ be an infinite finitely generated field, $M/K$ a separable
extension, PAC over $K$, $n$ a positive integer, and $N$ a Galois
extension of $M$ of degree at most $n$ with Galois group $A$. Then
there exist finite extensions $K'\cont L$ of $K$ such that with
$M'=K'M$ and $N'=K'N$ the following hold:
\renewcommand{\labelenumi}{\textnormal{(\alph{enumi})}}
\begin{enumerate}
\item
$N'=LM'$ and $\Gal(N'/M')\isom A$. 
\item
$L/K'$ is a Galois extension and $\Gal(L/K')\isom S_n$.
\end{enumerate}
$$
\xymatrix{
 L \ar@{-}[rr] \ar@{-}[d]_{S_n}
&& N' \ar@{-}[dl]_A \ar@{-}[dd] \cr K' \ar@{-}[r] \ar@{-}[d] & M'
\ar@{-}[d] \cr K \ar@{-}[r] & \ar@{-}[r] M & N }
$$
\end{lemma}

\begin{proof} We break the proof into three parts.

\noindent\textsc{Part A:} \textsl{Transcendental extensions.} First
we apply Lemma~\ref{lem:1} to construct Diagram (\ref{eq1}). Then we
choose $x\in F_0$ integral over $K[t]$ with $F_0=K(t,x)$ and we let
$g\in K[T,X]$ be the monic polynomial in $X$ such that
$g(t,X)=\irr(x,K(t))$. In particular, $r_1(t)=\discr(g(t,X))\in
K[t]$ and $r_1(t)\ne0$. Finally we choose $z\in D$ integral over
$M[t]$ with $D=M(t,z)$ and we let $f\in K[T,X]$ be the monic
polynomial such that $f(t,X)=\irr(z,M(t))$. Then
$r_2(t)=\discr(f(t,X))\in M[t]$ and $r_2(t)\ne0$. Since $D$ is
regular over $M$, the polynomial $f(T,X)$ is absolutely irreducible
\cite[Cor.~10.2.2(b)]{FrJ}. Let $r(t)=r_1(t)r_2(t)$.

Replacing $K$ by a finite extension in $M$, we may assume that $K$
contains all of the coefficients of $f(t,X)$, $g(t,X)$, and $r(t)$.
Set $R_0=K[t,r(t)^{-1}]$, $R=R_0M=M[t,r(t)^{-1}]$, $S_0=R_0[x]$,
$S=S_0M=R[x]$, and $V=R[z]$. Then $S_0/R_0$, $S/R$, and $V/R$ are
ring covers and $F_0/K(t)$, $F/M(t)$, and $D/M(t)$ are the
corresponding field covers.
$$
\xymatrix@=10pt{ S_0 \ar@{-}[rr] \ar@{-}[dd] && S \ar@{-}[dd]
\\
&&& V \ar@{-}[dl]
\\
R_0 \ar@{-}[rr] \ar@{-}[dd] && R \ar@{-}[dd]
\\
\\
K \ar@{-}[rr] && M \ar@{-}[rr] && N }
$$

\noindent\textsc{Part B:} \textsl{Specialization.} Lemma~\ref{lem:2}
gives a finite purely inseparable extension $K'$ of $K$, a
nonconstant rational function $q\in K'(T)$, and a finite subset $B$
of $K'$ such that $f(q(T),X)$ is absolutely irreducible, $g(q(a),X)$
is irreducible in $K'[X]$, and $r(q(a))\ne0$ for each $a\in
K'\hefresh B$.

We put $'$ on rings and fields to denote their composition with
$K'$. For example $M'=K'M$. Since $\Fhat/K$ is separable and $K'/K$
is purely inseparable, these extensions are linearly disjoint. It
follows that (a), (b), and (c) of Lemma~\ref{lem:1} hold for the
tagged rings and fields. In particular, $S'_0/R'_0$ and $S'/R'$ are
Galois covers of rings and $\Gal(N'/M')\isom A$. By
\cite[Cor.~2.5]{JaR1}, $M'$ is PAC over $K'$, hence there exists
$(a,c)\in K'\times M'$ such that $a\notin B$ and $f(b,c)=0$ with
$b=q(a)$. By the choice of $B$, \ $g(b,X)$ is irreducible in $K'[X]$
and $r(b)\ne0$.

The tag notation also gives $R'_0=K'[t,r(t)^{-1}]$ and
$V'=M'[t,z,r(t)^{-1}]$. Since $V'$ is integral over $R'$ and
$f(t,z)=0$, we may extend the specialization $(t,z)\to(b,c)$ to an
$M'$-epimorphism $\psi\colon V'\to M'$ satisfying $\psi(R'_0)=K'$.

\noindent\textsc{Part C:} \textsl{Finite extensions of $K'$.}  Let
$\Shat'$ be the integral closure of $R'$ in $\Fhat'$. Then
$S'=\Shat'\cap F'$. By \cite[Lemma 2.5.10]{FrJ},
$\Shat'=V'\tensor_{M'}N'$. Furthermore, $D'$ is linearly disjoint
from $F'$ over $F'\cap D'$, hence by the same lemma, $\Shat'=S'V'$.
$$
\xymatrix{ S'_0=R'_0[x] \ar@{-}[r] \ar@{-}[d] & S' \ar@{-}[r]
\ar@{-}[d] & \Shat' \ar@{-}[dl]|{V'=R'[z]} \ar@{-}[d]
\\
R'_0=K'[t,r(t)^{-1}] \ar@{-}[r] \ar@{-}[d] & R'=M'[t,r(t)^{-1}]
\ar@{-}[r] \ar@{-}[d] & R'N' \ar@{-}[d]
\\
K' \ar@{-}[r] & M' \ar@{-}[r] & N' } \qquad
$$
Setting $\psi(vn)=\psi(v)n$ for each $v'\in V'$ and $n\in N'$
extends $\psi$ to an $N'$-epimorphism $\psi\colon\Shat'\to N'$. In
particular, $M'\cont\psi(S')$, hence
$N'=\psi(\Shat')=\psi(S'V')=\psi(S')M'=\psi(S')$.

Let $L=K'(\psi(x))$. Then $\psi(S'_0)=\psi(R'_0[x])=K'(\psi(x))=L$.
Since $\psi(x)$ is a root of $g(b,X)$ and $g(b,X)$ is irreducible
over $K'$, we have
$$
[L:K']=\deg(g(b,X)) =\deg(g(t,X)) =[F_0:K(t)] =n!.
$$
Hence, $\Gal(L/K')\isom\Gal(F_0/K(t))\isom S_n$. Finally,
$N'=\psi(S')=\psi(S'_0M')=\psi(S'_0)M'=LM'$, as desired.
\end{proof}

We say that a separable algebraic extension $M$ of a field $K$ is
{\bf accessible} if there exists a sequence of fields
$$
K=K_0\cont K_1\cont K_2\cont\cdots\cont M
$$
such that $K_{i+1}/K_i$ is Galois for each $i$ and
$\union_{i=0}^\infty K_i=M$. In particular, every Galois extension
of $K$ is accessible. If $L/K$ is a finite Galois extension, then
the sequence $\Gal(L/L\cap K_i)$, $i=0,1,2,\ldots$, of subgroups of
$\Gal(L/K)$ is finite, so there is a positive integer $m$ such that
\begin{eqnarray*}
\lefteqn{\Gal(L/L\cap M) =}\\
&=&\Gal(L/L\cap K_m) \normal\Gal(L/L\cap K_{m-1})
\normal\cdots\normal\Gal(L/L\cap K_1) \normal\Gal(L/K).
\end{eqnarray*}
In other words, $\Gal(L/L\cap M)$ is a {\bf subnormal subgroup} of
$\Gal(L/K)$.

\begin{theorem}\label{thm:4}
Let $K$ be a finitely generated field, $M_0$ an accessible
extension of $K$, and $M$ a separable algebraic extension of $M_0$.
If $M$ is PAC over $K$ and $M\ne K_s$, then, as a supernatural
number, $[M:M_0]=\prod_pp^\infty$.
\end{theorem}

\begin{proof} By [JaR1, Remark 1.2(b)], $K$ is an infinite field.
Choose a proper finite Galois extension $N$ of $M$ with Galois group
$A$ and let $n$ be a positive integer dividing $5|A|$. Let $K'$ and
$L$ be fields satisfying Conditions (a) and (b) of
Lemma~\ref{lem:3}. Set $M_0'=K'M_0$, $M'=K'M$,  $L_0=L\cap M'_0$,
and $L_1 = L\cap M'$. Then $M_0'$ is an accessible extension of
$K'$, hence $\Gal(L/L_0)$ is a subnormal subgroup of
$\Gal(L/K')\isom S_n$. Since $n\geq 5$, the sequence $1\normal A_n
\normal S_n$ is the only composition series of $S_n$ \cite[p. 173,
Thm. 5.1]{Hup}. Therefore $\Gal(L/L_0)$ is either $1$ or $A_n$ or
$S_n$. By Condition (a) of Lemma~\ref{lem:3},
$$
\Gal(L/L_0) \ge \Gal(L/L_1) \cong \Gal(N'/M')\cong A \neq 1.
$$
Therefore, $A_n\le \Gal(L/L_0)$, so $\frac{(n-1)!}{2}$ divides
$[L_1:L_0]=[L_1M_0':M_0']$ (note that $\frac{|A_n|}{|A|}\mid
\frac{(n-1)!}{2}$).
$$
\xymatrix{ & L \ar@{-}[rrr] \ar@{-}[d]_A &&& N' \ar@{-}[ddd]
\ar@{-}[dl]_A
\\
& L_1 \ar@{-}[r] \ar@{-}[d] & L_1M'_0 \ar@{-}[d] \ar@{-}[r] & M'
\ar@{-}[dd]
\\
K' \ar@{-}[r] \ar@{-}[d] & L_0 \ar@{-}[r] & M'_0 \ar@{-}[d]
\\
K \ar@{-}[rr] && M_0 \ar@{-}[r] & M \ar@{-}[r] & N }
$$
Since $n$ is arbitrarily large and $L_1M_0'\cont M'$, we have
$[M':M_0']=\prod_{p}p^{\infty}$. Since $K'$ is a finite extension,
$[M:M_0]=\prod_{p}p^{\infty}$.
\end{proof}

The main result of this note is a special case of
Theorem~\ref{thm:4}:

\begin{corollary}
Let $K$ be a finitely generated field and let $N/K$ be a separable
extension PAC over $K$, \ $N\ne K_s$. Then $N$ is not an accessible
extension of $K$. In particular $N$ is not Galois over $K$.
\end{corollary}

\begin{corollary}
Let $K$ be an infinite finitely generated field and let $e$ be a
positive integer. Then, for almost all $\bfsig\in\Gal(K)^e$ the
extension $K_s(\bfsig)$ of $K$ is inaccessible.
\end{corollary}

\begin{proof} By \cite[Prop.~3.1]{JaR1}, for almost all
$\bfsig\in\Gal(K)^e$ the field $K_s(\bfsig)$ is PAC over $K$; in
addition, $K_s(\bfsig)\ne K_s$. Hence, by Theorem~\ref{thm:4},
$K_s(\bfsig)$ is inaccessible over $K$. 
\end{proof}

\begin{conjecture}
Let $K$ be a finitely generated field and $M$ an algebraic extension
of $K$. If $M$ is PAC over $K$, then $\Gal(M)$ is finitely
generated.
\end{conjecture}

%
%

\end{document}